\documentclass[letterpaper,fleqn,reqno,12pt]{amsart}
\usepackage[letterpaper,hmargin=1.25in,vmargin=1in]{geometry}
\usepackage{natbib}
\usepackage{rotating}

\newcommand{\CC}{\mathbb{C}}
\newcommand{\QQ}{\mathbb{Q}}
\newcommand{\ZZ}{\mathbb{Z}}
\newcommand{\set}[1]{\left\{#1\right\}}

\newtheorem{theorem}{Theorem}[section]

\begin{document}

\title[E. Carlini, G. Pistone]{Hilbert Bases for Orthogonal Arrays}

\author[E. Carlini, G. Pistone]{Enrico Carlini}
\author[]{Giovanni Pistone}

\address{DIMAT Politecnico di Torino}
\address{DIMAT Politecnico di Torino}

\email{enrico.carlini@polito.it}
\email{giovanni.pistone@polito.it}


\date{\today}


\begin{abstract}
In this paper, we relate the problem of generating all 2-level orthogonal arrays of given dimension and force, i.e. elements in OA$(n,m)$, where $n$ is the number of factors and $m$ the force, to the solution of an Integer Programming problem involving rational convex cones. We do not restrict the number of points in the array, i.e. we admit any number of replications. This problem can be theoretically solved by means of Hilbert bases which form a finite generating set for all the elements in in the infinite set OA$(n,m)$. We discuss some examples which are explicitly solved with a software performing Hilbert bases computation.
\end{abstract}

\maketitle
\section{Introduction}

We shall investigate orthogonal arrays with 2 levels coded $-1$, $+1$. In this paper, a full factorial design $D(n)$ is
the Cartesian product $\{-1,+1\}^n$ and a multi-subset $F$ of $D(n)$ is called a fraction. If $D(n)$ is embedded in the affine
space  $\QQ^n$, it is the set of solutions of the system of polynomial
equations $x_1^2 - 1 = 0$, $i = 1,\dots,n$. In a series of paper,
starting with \cite{pistone|wynn:1996}, this approach has been
systematically developed, mainly because of the availability of
algorithms and software able to solve systems of polynomial equations
with rational coefficients. The state of the art on year 2000 in
discussed in the book \cite{pistone|riccomagno|wynn:2001}.

In modern geometric terms, a fraction is a non-reduced 0-dimensional \emph{scheme} supported
on a subset of $D(n)$. The geometry of such schemes has probably many interesting thing to say to statistical design theory, but this will be the subject of further work. Here, we take a simpler approach, consisting in the remark that a fraction with replicates is  fully described by a functions $R : D(n)$ giving the number of replications $R(a) = 0,1,2,\dots$ of each point $a \in D(n)$. If the points of the full factorial design are listed in some order, then $R$ is a vector with non negative integer elements. A fraction of $D(n)$ such that all of its
$m$-dimensional orthogonal projections are replications of a full factorial designs
of the same cardinality (as multi-sets) is called an
$n$-dimensional orthogonal array of force $m$. We denote with
OA$(n,m)$ the set of such fractions.

One possible algebraic approach to the problem of describing the elements of
OA$(n,m)$ uses polynomial counting functions, i.e. polynomial
functions that take non negative integer values on $D(n)$. A
polynomial counting function on $D(n)$ that takes only the value 0 and
1 is called an indicator polynomial function. An indicator polynomial
function is associated with a fraction with no replications. This approach was developed in \cite{fontana|pistone|rogantin:96}, \cite{fontana|pistone|rogantin:2000}, \cite{ye:2003} for 2-level design and generalized in \cite{ye:2004}, \cite{pistone|rogantin:2006-17}.

In order to implement this kind of study, we need to describe the
ring of $\QQ$-valued functions over $D(n)$. This can be done by
means of standard techniques in Algebraic Statistic, see e.g
\cite{pistone|riccomagno|wynn:2001}. Let $\QQ$ be the rational
number field and consider the affine space $\QQ^n$. By realizing
$D(n)$ as a subset of $\QQ^n$, we can describe it using algebraic
equations. More precisely, $D(n)$ is associated to the
\emph{ideal} $I\subset R=\QQ[x_1,\ldots,x_n]$ generated by
\[
x_i^2-1=0, \quad i=1,\ldots,n.
\]
The ring of function over $D(n)$ is then the quotient ring $R/I$. As a $\QQ$-vector space,
$R/I$ is generated by the following set of monomials:
\begin{equation} \label{eq:monbasis}
\set{X^\alpha : \alpha=(\alpha_1,\ldots,\alpha_n) \in L = \set{0,1}^2}
\quad .
\end{equation}
The capital $X$ is intended to distinguish between the indeterminate $x$ and the corresponding function defined on $D(n)$, i.e. $X^\alpha(a) = \prod_{i=1}^n a_i^{\alpha_i}$, $a \in D(n)$.

Given a fraction $F$ of $D(n)$, its counting function can be
written as
\[R_F=\sum_\alpha b_\alpha X^\alpha\]
for a unique choice of constants $b_\alpha$'s. Notice that $R_F$
uniquely determines the fraction.

The previous setting is fully general and applies to any finite subset
of the affine space $\QQ^n$. However, the case we are discussing is
quite special, essentially because the mapping $\alpha \mapsto
\prod_{i=1}^n (-1)^{\alpha_i}$ is a $\CC$-representation of $L$ as the
additive group $\mod 2$. Then, the monomial functions in the vector basis in \eqref{eq:monbasis} are \emph{orthogonal}:
\[ \sum_{a \in D(n)} X^\alpha(a) X^\beta(a) =
\begin{cases}
0 & \text{ if $\alpha \ne \beta$, and} \\ 2^n & \text{ if $\alpha=\beta$.}
\end{cases}
\]
In particular, $\sum_{a \in D(n)} X^\alpha(a) = 0$, $\alpha \in L_0 = L \setminus \set{0}$. Then, the projection of the fraction $F$ on the factors with indices in $J \subset \set{1,\dots,n}$ is defined on $D(J)$ with counting function
\[
R_J(a_j) = \sum_{b_J \in D(J^\mathrm{c})} \sum_{\alpha = \alpha_J + \beta_J \in L = L_J \times L_{J^{\mathrm{c}}}} X^{\alpha_J}(a_J) X^{\beta_J}(b_J) = \sum_{\alpha_J \in L_j} X(a_j)^{\alpha_J}
\]
The conclusion is that condition $F\in \text{OA}(n,m)$ can be
easily expressed in terms of the $b_\alpha$'s. To see this, fix a
subset of indexes $J \subset \set{1,\ldots,n}$ of cardinality $m$. As $R_J$ has to be constant to be the counting function of a replication of a full factorial,then the orthogonality condition is expressed imposing the
vanishing of all $b$'s other than $b_0$ in this expression for each such $J$. Namely we get
conditions: \emph{$b_\alpha=0$, $\alpha \ne 0$ has no more than $m$ non-zero
component}.

A counting polynomial $R$ is a polynomial which has non negative integer values on $D(n)$. An algebraic way to say that is the following. Let $f_n(x) = x (x-1) \cdots (x - n)$ be the factorial polynomial of order $n$. A polynomial $R$ takes values in $\ZZ_+$ if all but a finite number of the polynomials $f(R_n(x))$ are identically zero on $D(n)$. In particular, the values are 0 and 1 if and only if $R(R-1) = 0$ on $D(n)$, or $\sum_{\alpha+\beta=\gamma} b_\alpha b_\beta = 0$ for all $\gamma \in L$. Thus we have a method to find the elements of OA$(n,m)$ without replications,
namely we have to solve in the $b_\alpha$'s the system
\[
\begin{cases}
  \sum_{\alpha+\beta=\gamma} b_\alpha b_\beta = 0 & \gamma \in L \\ b_\alpha = 0 & 1 \leq \sum \alpha \leq m
\end{cases}
\]
This was done in \cite{fontana|pistone|rogantin:2000} using the software \texttt{CoCoA}. For future reference we reproduce here the results for the case OA$(5,2)$. This method founds 1054 fractions, that were classified as follows. Below, $C_I$ denote the term $X^\alpha$ with $\alpha(i) = 1$ for $i \in I$ and zero otherwise.
\begin{enumerate}
\item
92 are regular fractions;  there exist 3 classes of equivalence for change of
signs and permutation of factors. Among those,
\begin{enumerate}
\item
32 have $b_0=1/2$, i.e. 16 points:
\begin{enumerate}
\item
2 have  counting polynomial of the form:
$$\frac 1 2 + \frac 1 2 C_{ijhkl} \quad ;
$$
and the fraction has resolution 5 and projectivity 4.
\item
10 have counting polynomial of the form:
$$\frac 1 2 + \frac 1 2 C_{ijhk} \quad ;
$$
the fraction has resolution 4 and projectivity 3; moreover it
fully projects on the $\{jhkl\}$-factors, on the $\{ihkl\}$-factors, on the
$\{ijkl\}$-factors and on the $\{ijhl\}$-factors.
\item
20 have counting polynomial of the form:
$$\frac 1 2 + \frac 1 2 C_{ijh} \quad ;
$$
the corresponding fraction has resolution 3 and projectivity 2; moreover it
fully projects on the $\{jhkl\}$-factors, on the $\{ihkl\}$-factors and on the
$\{ijkl\}$-factors.
\end{enumerate}
\item
 60 have $b_0=1/4$ (8 points) and resolution 3.  This is an unique class of equivalence
for change of signs and permutation of factors.
A counting polynomial is of the form:
$$\frac 1 4+ \frac 1 4 C_{ijh}+ \frac 1 4 C_{hkl}+ \frac 1 4 C_{ijkl} \quad .
$$
\end{enumerate}
\item
  940 are non regular fractions with resolution 3:
\begin{enumerate}
\item
192 have $b_0=3/8$ (12 points); this is an unique class of  equivalence for
change of sign and permutation of factors. The counting polynomials have all the
coefficients of the interactions of order 3 and 4 equal, in absolute value,
to $\frac 1 8$ and the coefficient of the interaction of order 5 equal to 0.
\item
520 have $b_0=1/2$ (16 points); there exist 4 classes of equivalence for change of
sign and permutation of factors:
\begin{enumerate}
\item
a counting polynomial is:
$$\frac 1 2 + \frac 1 4 C_{ijh}+ \frac 1 4 C_{ijhk}+ \frac 1 4 C_{ijhl}-\frac
1 4 C_{ijhkl} \quad ;
$$
the corresponding fraction fully projects on the ${jhkl}$-factors, on the
${ihkl}$-factors and on the ${ijkl}$-factors.
\item
a counting polynomial is:
$$\frac 1 2 + \frac 1 4 C_{ijh}+ \frac 1 4 C_{ijk}+ \frac 1 4 C_{ijhl}-\frac
1 4 C_{ijkl} \quad ;
$$
the fraction fully projects on the ${jhkl}$-factors  and on the
${ihkl}$-factors.
\item
an counting polynomial is:
$$\frac 1 2 + \frac 1 4 C_{ijh}+ \frac 1 4 C_{ijk}+ \frac 1 4 C_{ijl}-\frac
1 4 C_{ijhkl} \quad ;
$$
the fraction fully projects on the ${jhkl}$-factors  and on the
${ihkl}$-factors.
\item
a counting polynomial is:
$$\frac 1 2 + \frac 1 4 C_{ijh}+ \frac 1 4 C_{ijk}+ \frac 1 4 C_{ihl}-\frac
1 4 C_{ikl} \quad ;
$$
the fraction fully projects on the ${jhkl}$-factors.
\end{enumerate}
\item
192 have $b_0=5/8$; they are the complement of the fractions with  $b_0=3/8$.
\item
60 have $b_0=3/4$; they are the complement of the fractions with  $b_0=1/4$.
\end{enumerate}
\end{enumerate}
See the original paper for more details. Discarding the 192 + 60 = 252 fractions that are complement of a smaller one, in conclusion there are 1054 - 252 = 802 fractions, 60 with 8 points, 192 with 12 points, 552 with 16 points.
\section{Cones and semi-groups}
We follow the presentation of \cite{schrijver:1986}. Let $C$ be a (rational)
\emph{convex cone} in $\QQ^n$, i.e. $C \subset \QQ^n$ and $x,y \in
C$, $\lambda, \mu \in \QQ_+$ imply $\lambda x + \mu y \in C$. A convex
cone is a convex subset of $\QQ^n$, so we say simply \emph{cone}. The
cone $C$ is \emph{polyhedral} if the following two equivalent
conditions are satisfied:
\begin{enumerate}
\item \label{item:polycone} $C = \set{x \in \QQ : Ax \geq 0}$, for
  some matrix $A \in \QQ^{k,n}$;
\item \label{item:fingen} $C = \set{\lambda_1 x_1 + \cdots + \lambda_m x_m :
    \lambda_1,\dots,\lambda_m \in \QQ_+}$, for some finite set of
  vectors $x_1,\dots,x_m \in \QQ^m$, called \emph{generating set} of $C$.
\end{enumerate}
Notice that the direction of the inequality in Item
\ref{item:polycone} is unessential because we can change $A$ to $-A$,
and that we could have equality, represented with a $2k\times m$
matrix with blocks $A$ and $-A$. The cone is said to be
\emph{pointed} if $C \cap -C = \set{0}$, which in turn is equivalent
to $Ay=0$ implies $y = 0$

Let $O = C \cap \ZZ^n$ be the set of lattice points of the cone
$C$. Then $O$ is a sub-semi-group of $\ZZ^n$ because $r,s \in O$
implies $r + s \in O$. If $x_1,\dots,x_m$ is a set of generators of
the cone $C$, than each lattice point $r \in O$ can be written as $r =
\lambda_1 x_1 + \cdots + \lambda_m x_m$ with rational coefficients.

A \emph{Hilbert basis} of $O = C \cap \ZZ^n$ is a finite set of
elements $r_1,\dots,r_l$ such that any other element of $O$ is
linear combination with non-negative integer coefficients of the
$r_i$'s, see \cite[Sec. 16.4]{schrijver:1986}. In other words, a
Hilbert basis is a generating set of $O$ as a semi-group.

\begin{theorem}[Existence and uniqueness of Hilbert bases] Each
  rational polyhedral cone is generated by an Hilbert basis of its
  lattice points. If $C$ is
  pointed, then there exist a unique inclusion-minimal  Hilbert basis. \end{theorem}
The previous Theorem does not apply to more general sub-semi-groups of
$\ZZ^n$. In fact, the following was proved by Hemmecke and
Weissmantel.
\begin{theorem}
  A sub-semi-group $S$ of $\ZZ^n$ has a finite generating set if and
  only if the cone of $S$ is polyhedral.
\end{theorem}
We are going to use the previous theory in a special case. We
are given a rational matrix $M_1 \in \QQ^{N,k}$ and the cone of
non-negative $y \in \QQ_+^N$ such that $y^t M_1 = 0$. This is a
polyhedral cone defined by $M_1^t y = 0$ and $Iy \geq 0$. Notice
that the vector space $\ker M_1^t$ has a linear basis of integral
vectors. Moreover, if $1^t M_1 = 0$, then there exist a linear
basis of non-negative integral vectors. This basis allows us to
obtain all lattice points using rational coefficients while the
Hilbert basis, being bigger, does the same using non-negative
integer coefficients.
\section{Fractions and semi-groups}

An approach to the generation of all OA$(n,m)$ involves Hilbert bases
of a semi-group. Let $R(a)$, $a \in D(n)$ be a vector of replicates,
$R(a) \in \ZZ_+$. As an integer valued polynomial of the ring
$\QQ[x_1,\dots,x_n]$, reduced on the ideal $I(n)$, $R$ can be written
on $D(n)$ uniquely as
\[R(x) = \sum_{\alpha \in L} b_\alpha x^\alpha,\] $L = \set{0,1}^n$.
Notice that, for a fixed $\alpha\in L$, we have
\[
\sum_{a\in D(n)} R(a)X^\alpha(a)=2^n\cdot b_\alpha.
\]
Thus we have a linear system of equations in the vector of replicates
$R=(\ldots,F(a),\ldots)$ as $\alpha$ varies in $L$. Assume now we have
a subset $L_1 \subset L_0 = L \setminus 0$, and we want $b_\alpha = 0$
for all $\alpha \in L_1$. Then we are looking for the non negative
integer solutions of the system
\begin{equation}
  \label{eq:system}
  \sum_{a\in D(n)} R(a)X^\alpha(a)=0 \quad \alpha \in L_1
\end{equation}
The model matrix restricted to $L_1$, namely
$M_1=(\mathbf{x}^\alpha(a))_{a\in D(n),\alpha\in L_1}$, is the matrix
of coefficients in Equation \eqref{eq:system}. The non negative integer
solutions of the integer linear system of equations $R^t\cdot M=0$ are
the lattice points of the cone $C = \ker M_1^t$. Hence we look for
generators of the semi-group $O =C \cap \ZZ_+^n$.

A \emph{Hilbert basis} of $O$ is a finite and minimal set of
elements $r_1,\ldots,r_l$ such that any other element of $C_+$ is
linear combination with non-negative integer coefficients of the
$r_i$'s. The list $L_1$ consists of all interactions of order
between 1 and $m$, $r_i$ gives us an element in OA$(n,m)$ and the
integers combination $\sum n_i r_i$, the fraction obtained taking
the union of the corresponding fractions $n_i$ times each, is the
generic element of OA$(n,m)$. Notice that the matrix $M_1^t$ has
entries equal to either $1$ or $-1$, hence $\ker M_1^t$, as a
$\QQ$-vector space has a basis of $2^n - \# L_1$ non negative
integer vectors.  These basis elements produce elements in
OA$(n,m)$ which are independent, i.e. one does not decompose as
union of the others. However, in general this is not a Hilbert
basis as these fractions can be decomposed as union of other
elements in OA$(n,m)$.

Algorithms for computing Hilbert bases are implemented in
specialized software, namely {\tt 4ti2}, see
\cite{hemmecke|hemmecke|malkin:2005} and {\tt CoCoA}, see
\cite{CocoaSystem}. Using the former we present some examples of
interest in order to discuss the use of Hilbert bases in designs
theory. We remark that the computation of Hilbert bases with {\tt
4ti2} uses the Project-and-Lift algorithm described in
\cite{hemmecke:2006}.

\section{Examples}

\subsection{OA$(5,3)$}
In this case the $32\times 25$ matrix $M_1$ is given in Table
\ref{tab:1}%
\begin{table}
  \centering
  \caption{Matrix $M_1$ for the OA$(5,3)$}
  \label{tab:1}
  \tiny
  \begin{sideways}
  $\begin{array}{ rrrrrrrrrrrrrrrrrrrrrrrrrr }%
----- & -1 & -1 &  1 & -1 &  1 &  1 & -1 & -1 &  1 &  1 & -1 &  1 & -1 & -1 & -1 &  1 &  1 & -1 &  1 & -1 & -1 &  1 & -1 & -1 & -1\\
----+ &  1 & -1 & -1 & -1 & -1 &  1 &  1 & -1 & -1 &  1 &  1 &  1 &  1 & -1 & -1 & -1 &  1 &  1 &  1 &  1 & -1 &  1 &  1 & -1 & -1\\
---+- & -1 &  1 & -1 & -1 &  1 & -1 &  1 & -1 &  1 & -1 &  1 &  1 & -1 &  1 & -1 &  1 & -1 &  1 &  1 & -1 &  1 &  1 & -1 &  1 & -1\\
---++ &  1 &  1 &  1 & -1 & -1 & -1 & -1 & -1 & -1 & -1 & -1 &  1 &  1 &  1 & -1 & -1 & -1 & -1 &  1 &  1 &  1 &  1 &  1 &  1 & -1\\
--+-- & -1 & -1 &  1 &  1 & -1 & -1 &  1 & -1 &  1 &  1 & -1 & -1 &  1 &  1 & -1 &  1 &  1 & -1 & -1 &  1 &  1 &  1 & -1 & -1 &  1\\
--+-+ &  1 & -1 & -1 &  1 &  1 & -1 & -1 & -1 & -1 &  1 &  1 & -1 & -1 &  1 & -1 & -1 &  1 &  1 & -1 & -1 &  1 &  1 &  1 & -1 &  1\\
--++- & -1 &  1 & -1 &  1 & -1 &  1 & -1 & -1 &  1 & -1 &  1 & -1 &  1 & -1 & -1 &  1 & -1 &  1 & -1 &  1 & -1 &  1 & -1 &  1 &  1\\
--+++ &  1 &  1 &  1 &  1 &  1 &  1 &  1 & -1 & -1 & -1 & -1 & -1 & -1 & -1 & -1 & -1 & -1 & -1 & -1 & -1 & -1 &  1 &  1 &  1 &  1\\
-+--- & -1 & -1 &  1 & -1 &  1 &  1 & -1 &  1 & -1 & -1 &  1 & -1 &  1 &  1 & -1 &  1 &  1 & -1 &  1 & -1 & -1 & -1 &  1 &  1 &  1\\
-+--+ &  1 & -1 & -1 & -1 & -1 &  1 &  1 &  1 &  1 & -1 & -1 & -1 & -1 &  1 & -1 & -1 &  1 &  1 &  1 &  1 & -1 & -1 & -1 &  1 &  1\\
-+-+- & -1 &  1 & -1 & -1 &  1 & -1 &  1 &  1 & -1 &  1 & -1 & -1 &  1 & -1 & -1 &  1 & -1 &  1 &  1 & -1 &  1 & -1 &  1 & -1 &  1\\
-+-++ &  1 &  1 &  1 & -1 & -1 & -1 & -1 &  1 &  1 &  1 &  1 & -1 & -1 & -1 & -1 & -1 & -1 & -1 &  1 &  1 &  1 & -1 & -1 & -1 &  1\\
-++-- & -1 & -1 &  1 &  1 & -1 & -1 &  1 &  1 & -1 & -1 &  1 &  1 & -1 & -1 & -1 &  1 &  1 & -1 & -1 &  1 &  1 & -1 &  1 &  1 & -1\\
-++-+ &  1 & -1 & -1 &  1 &  1 & -1 & -1 &  1 &  1 & -1 & -1 &  1 &  1 & -1 & -1 & -1 &  1 &  1 & -1 & -1 &  1 & -1 & -1 &  1 & -1\\
-+++- & -1 &  1 & -1 &  1 & -1 &  1 & -1 &  1 & -1 &  1 & -1 &  1 & -1 &  1 & -1 &  1 & -1 &  1 & -1 &  1 & -1 & -1 &  1 & -1 & -1\\
-++++ &  1 &  1 &  1 &  1 &  1 &  1 &  1 &  1 &  1 &  1 &  1 &  1 &  1 &  1 & -1 & -1 & -1 & -1 & -1 & -1 & -1 & -1 & -1 & -1 & -1\\
+---- & -1 & -1 &  1 & -1 &  1 &  1 & -1 & -1 &  1 &  1 & -1 &  1 & -1 & -1 &  1 & -1 & -1 &  1 & -1 &  1 &  1 & -1 &  1 &  1 &  1\\
+---+ &  1 & -1 & -1 & -1 & -1 &  1 &  1 & -1 & -1 &  1 &  1 &  1 &  1 & -1 &  1 &  1 & -1 & -1 & -1 & -1 &  1 & -1 & -1 &  1 &  1\\
+--+- & -1 &  1 & -1 & -1 &  1 & -1 &  1 & -1 &  1 & -1 &  1 &  1 & -1 &  1 &  1 & -1 &  1 & -1 & -1 &  1 & -1 & -1 &  1 & -1 &  1\\
+--++ &  1 &  1 &  1 & -1 & -1 & -1 & -1 & -1 & -1 & -1 & -1 &  1 &  1 &  1 &  1 &  1 &  1 &  1 & -1 & -1 & -1 & -1 & -1 & -1 &  1\\
+-+-- & -1 & -1 &  1 &  1 & -1 & -1 &  1 & -1 &  1 &  1 & -1 & -1 &  1 &  1 &  1 & -1 & -1 &  1 &  1 & -1 & -1 & -1 &  1 &  1 & -1\\
+-+-+ &  1 & -1 & -1 &  1 &  1 & -1 & -1 & -1 & -1 &  1 &  1 & -1 & -1 &  1 &  1 &  1 & -1 & -1 &  1 &  1 & -1 & -1 & -1 &  1 & -1\\
+-++- & -1 &  1 & -1 &  1 & -1 &  1 & -1 & -1 &  1 & -1 &  1 & -1 &  1 & -1 &  1 & -1 &  1 & -1 &  1 & -1 &  1 & -1 &  1 & -1 & -1\\
+-+++ &  1 &  1 &  1 &  1 &  1 &  1 &  1 & -1 & -1 & -1 & -1 & -1 & -1 & -1 &  1 &  1 &  1 &  1 &  1 &  1 &  1 & -1 & -1 & -1 & -1\\
++--- & -1 & -1 &  1 & -1 &  1 &  1 & -1 &  1 & -1 & -1 &  1 & -1 &  1 &  1 &  1 & -1 & -1 &  1 & -1 &  1 &  1 &  1 & -1 & -1 & -1\\
++--+ &  1 & -1 & -1 & -1 & -1 &  1 &  1 &  1 &  1 & -1 & -1 & -1 & -1 &  1 &  1 &  1 & -1 & -1 & -1 & -1 &  1 &  1 &  1 & -1 & -1\\
++-+- & -1 &  1 & -1 & -1 &  1 & -1 &  1 &  1 & -1 &  1 & -1 & -1 &  1 & -1 &  1 & -1 &  1 & -1 & -1 &  1 & -1 &  1 & -1 &  1 & -1\\
++-++ &  1 &  1 &  1 & -1 & -1 & -1 & -1 &  1 &  1 &  1 &  1 & -1 & -1 & -1 &  1 &  1 &  1 &  1 & -1 & -1 & -1 &  1 &  1 &  1 & -1\\
+++-- & -1 & -1 &  1 &  1 & -1 & -1 &  1 &  1 & -1 & -1 &  1 &  1 & -1 & -1 &  1 & -1 & -1 &  1 &  1 & -1 & -1 &  1 & -1 & -1 &  1\\
+++-+ &  1 & -1 & -1 &  1 &  1 & -1 & -1 &  1 &  1 & -1 & -1 &  1 &  1 & -1 &  1 &  1 & -1 & -1 &  1 &  1 & -1 &  1 &  1 & -1 &  1\\
++++- & -1 &  1 & -1 &  1 & -1 &  1 & -1 &  1 & -1 &  1 & -1 &  1 & -1 &  1 &  1 & -1 &  1 & -1 &  1 & -1 &  1 &  1 & -1 &  1 &  1\\
+++++ &  1 &  1 &  1 &  1 &  1 &  1 &  1 &  1 &  1 &  1 &  1 &  1 &  1 &  1 &  1 &  1 &  1 &  1 &  1 &  1 &  1 &  1 &  1 &  1 &  1
\end{array}$
\end{sideways}
 \normalsize
 \end{table}
The matrix was generated with the software R using the function
computing the generalized Kronecker product of two arrays. Then, the
treatment points and the interactions are listed in right-to-left
lexicographic order. The $\alpha$'s of the interactions of order 1,2,3
are printed in the order: 00001, 00010, 00011, 00100, 00101, 00110,
00111, 01000, 01001, 01010, 01011, 01100, 01101, 01110, 10000, 10001,
10010, 10011, 10100, 10101, 10110,11000, 11001, 11010,
11100. Treatment values are shown; column names are not printed to
save space.

One run on my slow laptop of the \texttt{4ti2} function
\texttt{hilbert} has taken $0.05$ seconds to find the 28 elements in the
minimal Hilbert basis of OA$(5,3)$, see Table \ref{tab:2}.
\begin{table}
  \centering
  \caption{Hilbert basis for the OA$(5,3)$}
  \label{tab:2}
   \tiny
  \begin{sideways}
    $\begin{array}{ rrrrrrrrrrrrrrrrrrrrrrrrrrrrr }%
 & 1 & 2 & 3 & 4 & 5 & 6 & 7 & 8 & 9 & 10 & 11 & 12 & 13 & 14 & 15 & 16 & 17 & 18 & 19 & 20 & 21 & 22 & 23 & 24 & 25 & 26 & 27 & 28 \\
----- & 1 & 1 & 1 & 0 & 1 & 0 & 1 & 0 & 1 & 0 & 0 & 0 & 1 & 2 & 1 & 1 & 1 & 1 & 1 & 0 & 1 & 0 & 0 & 1 & 1 & 1 & 0 & 0\\
----+ & 1 & 0 & 0 & 1 & 0 & 1 & 0 & 1 & 0 & 1 & 1 & 0 & 1 & 0 & 1 & 1 & 1 & 1 & 0 & 2 & 0 & 1 & 1 & 1 & 0 & 0 & 1 & 1\\
---+- & 0 & 0 & 0 & 1 & 0 & 0 & 1 & 1 & 0 & 1 & 1 & 1 & 0 & 0 & 1 & 0 & 1 & 0 & 1 & 1 & 0 & 1 & 1 & 1 & 1 & 1 & 2 & 1\\
---++ & 0 & 1 & 1 & 0 & 1 & 1 & 0 & 0 & 1 & 0 & 0 & 1 & 1 & 1 & 0 & 1 & 0 & 1 & 1 & 0 & 2 & 1 & 1 & 0 & 1 & 1 & 0 & 1\\
--+-- & 0 & 0 & 0 & 0 & 0 & 1 & 0 & 1 & 1 & 1 & 1 & 1 & 0 & 0 & 0 & 1 & 0 & 1 & 0 & 1 & 1 & 1 & 1 & 1 & 1 & 1 & 1 & 2\\
--+-+ & 0 & 1 & 1 & 1 & 1 & 0 & 1 & 0 & 0 & 0 & 0 & 1 & 1 & 1 & 1 & 0 & 1 & 0 & 2 & 0 & 1 & 1 & 1 & 0 & 1 & 1 & 1 & 0\\
--++- & 1 & 1 & 1 & 1 & 1 & 1 & 0 & 0 & 0 & 0 & 0 & 0 & 2 & 1 & 1 & 1 & 1 & 1 & 1 & 1 & 1 & 1 & 1 & 0 & 0 & 0 & 0 & 0\\
--+++ & 1 & 0 & 0 & 0 & 0 & 0 & 1 & 1 & 1 & 1 & 1 & 0 & 0 & 1 & 1 & 1 & 1 & 1 & 0 & 1 & 0 & 0 & 0 & 2 & 1 & 1 & 1 & 1\\
-+--- & 0 & 1 & 0 & 1 & 0 & 1 & 0 & 1 & 0 & 1 & 0 & 1 & 1 & 0 & 1 & 1 & 0 & 0 & 1 & 1 & 1 & 2 & 1 & 0 & 1 & 0 & 1 & 1\\
-+--+ & 0 & 0 & 1 & 0 & 1 & 0 & 1 & 0 & 1 & 0 & 1 & 1 & 0 & 1 & 0 & 0 & 1 & 1 & 1 & 0 & 1 & 0 & 1 & 1 & 1 & 2 & 1 & 1\\
-+-+- & 1 & 0 & 1 & 0 & 1 & 1 & 0 & 0 & 1 & 0 & 1 & 0 & 1 & 1 & 0 & 1 & 1 & 2 & 0 & 1 & 1 & 0 & 1 & 1 & 0 & 1 & 0 & 1\\
-+-++ & 1 & 1 & 0 & 1 & 0 & 0 & 1 & 1 & 0 & 1 & 0 & 0 & 1 & 1 & 2 & 1 & 1 & 0 & 1 & 1 & 0 & 1 & 0 & 1 & 1 & 0 & 1 & 0\\
-++-- & 1 & 0 & 1 & 1 & 1 & 0 & 1 & 0 & 0 & 0 & 1 & 0 & 1 & 1 & 1 & 0 & 2 & 1 & 1 & 1 & 0 & 0 & 1 & 1 & 0 & 1 & 1 & 0\\
-++-+ & 1 & 1 & 0 & 0 & 0 & 1 & 0 & 1 & 1 & 1 & 0 & 0 & 1 & 1 & 1 & 2 & 0 & 1 & 0 & 1 & 1 & 1 & 0 & 1 & 1 & 0 & 0 & 1\\
-+++- & 0 & 1 & 0 & 0 & 0 & 0 & 1 & 1 & 1 & 1 & 0 & 1 & 0 & 1 & 1 & 1 & 0 & 0 & 1 & 0 & 1 & 1 & 0 & 1 & 2 & 1 & 1 & 1\\
-++++ & 0 & 0 & 1 & 1 & 1 & 1 & 0 & 0 & 0 & 0 & 1 & 1 & 1 & 0 & 0 & 0 & 1 & 1 & 1 & 1 & 1 & 1 & 2 & 0 & 0 & 1 & 1 & 1\\
+---- & 0 & 0 & 1 & 1 & 0 & 1 & 0 & 1 & 0 & 0 & 1 & 1 & 1 & 0 & 0 & 0 & 1 & 1 & 1 & 1 & 1 & 1 & 2 & 0 & 0 & 1 & 1 & 1\\
+---+ & 0 & 1 & 0 & 0 & 1 & 0 & 1 & 0 & 1 & 1 & 0 & 1 & 0 & 1 & 1 & 1 & 0 & 0 & 1 & 0 & 1 & 1 & 0 & 1 & 2 & 1 & 1 & 1\\
+--+- & 1 & 1 & 0 & 0 & 1 & 1 & 0 & 0 & 1 & 1 & 0 & 0 & 1 & 1 & 1 & 2 & 0 & 1 & 0 & 1 & 1 & 1 & 0 & 1 & 1 & 0 & 0 & 1\\
+--++ & 1 & 0 & 1 & 1 & 0 & 0 & 1 & 1 & 0 & 0 & 1 & 0 & 1 & 1 & 1 & 0 & 2 & 1 & 1 & 1 & 0 & 0 & 1 & 1 & 0 & 1 & 1 & 0\\
+-+-- & 1 & 1 & 0 & 1 & 1 & 0 & 1 & 0 & 0 & 1 & 0 & 0 & 1 & 1 & 2 & 1 & 1 & 0 & 1 & 1 & 0 & 1 & 0 & 1 & 1 & 0 & 1 & 0\\
+-+-+ & 1 & 0 & 1 & 0 & 0 & 1 & 0 & 1 & 1 & 0 & 1 & 0 & 1 & 1 & 0 & 1 & 1 & 2 & 0 & 1 & 1 & 0 & 1 & 1 & 0 & 1 & 0 & 1\\
+-++- & 0 & 0 & 1 & 0 & 0 & 0 & 1 & 1 & 1 & 0 & 1 & 1 & 0 & 1 & 0 & 0 & 1 & 1 & 1 & 0 & 1 & 0 & 1 & 1 & 1 & 2 & 1 & 1\\
+-+++ & 0 & 1 & 0 & 1 & 1 & 1 & 0 & 0 & 0 & 1 & 0 & 1 & 1 & 0 & 1 & 1 & 0 & 0 & 1 & 1 & 1 & 2 & 1 & 0 & 1 & 0 & 1 & 1\\
++--- & 1 & 0 & 0 & 0 & 1 & 0 & 1 & 0 & 1 & 1 & 1 & 0 & 0 & 1 & 1 & 1 & 1 & 1 & 0 & 1 & 0 & 0 & 0 & 2 & 1 & 1 & 1 & 1\\
++--+ & 1 & 1 & 1 & 1 & 0 & 1 & 0 & 1 & 0 & 0 & 0 & 0 & 2 & 1 & 1 & 1 & 1 & 1 & 1 & 1 & 1 & 1 & 1 & 0 & 0 & 0 & 0 & 0\\
++-+- & 0 & 1 & 1 & 1 & 0 & 0 & 1 & 1 & 0 & 0 & 0 & 1 & 1 & 1 & 1 & 0 & 1 & 0 & 2 & 0 & 1 & 1 & 1 & 0 & 1 & 1 & 1 & 0\\
++-++ & 0 & 0 & 0 & 0 & 1 & 1 & 0 & 0 & 1 & 1 & 1 & 1 & 0 & 0 & 0 & 1 & 0 & 1 & 0 & 1 & 1 & 1 & 1 & 1 & 1 & 1 & 1 & 2\\
+++-- & 0 & 1 & 1 & 0 & 0 & 1 & 0 & 1 & 1 & 0 & 0 & 1 & 1 & 1 & 0 & 1 & 0 & 1 & 1 & 0 & 2 & 1 & 1 & 0 & 1 & 1 & 0 & 1\\
+++-+ & 0 & 0 & 0 & 1 & 1 & 0 & 1 & 0 & 0 & 1 & 1 & 1 & 0 & 0 & 1 & 0 & 1 & 0 & 1 & 1 & 0 & 1 & 1 & 1 & 1 & 1 & 2 & 1\\
++++- & 1 & 0 & 0 & 1 & 1 & 1 & 0 & 0 & 0 & 1 & 1 & 0 & 1 & 0 & 1 & 1 & 1 & 1 & 0 & 2 & 0 & 1 & 1 & 1 & 0 & 0 & 1 & 1\\
+++++ & 1 & 1 & 1 & 0 & 0 & 0 & 1 & 1 & 1 & 0 & 0 & 0 & 1 & 2 & 1 & 1 & 1 & 1 & 1 & 0 & 1 & 0 & 0 & 1 & 1 & 1 & 0 & 0
\end{array}
$
  \end{sideways}\normalsize\end{table}
There are 12 arrays with 16 treatments an no replications; there are
16 orthogonal arrays with 24 points, all with replications.
\subsection{OA$(5,2)$}
In this case the $32 \times 15$ matrix $M_1$ is given in Table \ref{tab:3}.
\begin{table}
  \centering
  \caption{matrix $M_1$ for the OA$(5,2)$}
  \label{tab:3}
   \tiny \begin{sideways}$\begin{array}{ rrrrrrrrrrrrrrrr }%
 & 00001 & 00010 & 00011 & 00100 & 00101 & 00110 & 01000 & 01001 & 01010 & 01100 & 10000 & 10001 & 10010 & 10100 & 11000 \\
----- & -1 & -1 &  1 & -1 &  1 &  1 & -1 &  1 &  1 &  1 & -1 &  1 &  1 &  1 &  1\\
----+ &  1 & -1 & -1 & -1 & -1 &  1 & -1 & -1 &  1 &  1 & -1 & -1 &  1 &  1 &  1\\
---+- & -1 &  1 & -1 & -1 &  1 & -1 & -1 &  1 & -1 &  1 & -1 &  1 & -1 &  1 &  1\\
---++ &  1 &  1 &  1 & -1 & -1 & -1 & -1 & -1 & -1 &  1 & -1 & -1 & -1 &  1 &  1\\
--+-- & -1 & -1 &  1 &  1 & -1 & -1 & -1 &  1 &  1 & -1 & -1 &  1 &  1 & -1 &  1\\
--+-+ &  1 & -1 & -1 &  1 &  1 & -1 & -1 & -1 &  1 & -1 & -1 & -1 &  1 & -1 &  1\\
--++- & -1 &  1 & -1 &  1 & -1 &  1 & -1 &  1 & -1 & -1 & -1 &  1 & -1 & -1 &  1\\
--+++ &  1 &  1 &  1 &  1 &  1 &  1 & -1 & -1 & -1 & -1 & -1 & -1 & -1 & -1 &  1\\
-+--- & -1 & -1 &  1 & -1 &  1 &  1 &  1 & -1 & -1 & -1 & -1 &  1 &  1 &  1 & -1\\
-+--+ &  1 & -1 & -1 & -1 & -1 &  1 &  1 &  1 & -1 & -1 & -1 & -1 &  1 &  1 & -1\\
-+-+- & -1 &  1 & -1 & -1 &  1 & -1 &  1 & -1 &  1 & -1 & -1 &  1 & -1 &  1 & -1\\
-+-++ &  1 &  1 &  1 & -1 & -1 & -1 &  1 &  1 &  1 & -1 & -1 & -1 & -1 &  1 & -1\\
-++-- & -1 & -1 &  1 &  1 & -1 & -1 &  1 & -1 & -1 &  1 & -1 &  1 &  1 & -1 & -1\\
-++-+ &  1 & -1 & -1 &  1 &  1 & -1 &  1 &  1 & -1 &  1 & -1 & -1 &  1 & -1 & -1\\
-+++- & -1 &  1 & -1 &  1 & -1 &  1 &  1 & -1 &  1 &  1 & -1 &  1 & -1 & -1 & -1\\
-++++ &  1 &  1 &  1 &  1 &  1 &  1 &  1 &  1 &  1 &  1 & -1 & -1 & -1 & -1 & -1\\
+---- & -1 & -1 &  1 & -1 &  1 &  1 & -1 &  1 &  1 &  1 &  1 & -1 & -1 & -1 & -1\\
+---+ &  1 & -1 & -1 & -1 & -1 &  1 & -1 & -1 &  1 &  1 &  1 &  1 & -1 & -1 & -1\\
+--+- & -1 &  1 & -1 & -1 &  1 & -1 & -1 &  1 & -1 &  1 &  1 & -1 &  1 & -1 & -1\\
+--++ &  1 &  1 &  1 & -1 & -1 & -1 & -1 & -1 & -1 &  1 &  1 &  1 &  1 & -1 & -1\\
+-+-- & -1 & -1 &  1 &  1 & -1 & -1 & -1 &  1 &  1 & -1 &  1 & -1 & -1 &  1 & -1\\
+-+-+ &  1 & -1 & -1 &  1 &  1 & -1 & -1 & -1 &  1 & -1 &  1 &  1 & -1 &  1 & -1\\
+-++- & -1 &  1 & -1 &  1 & -1 &  1 & -1 &  1 & -1 & -1 &  1 & -1 &  1 &  1 & -1\\
+-+++ &  1 &  1 &  1 &  1 &  1 &  1 & -1 & -1 & -1 & -1 &  1 &  1 &  1 &  1 & -1\\
++--- & -1 & -1 &  1 & -1 &  1 &  1 &  1 & -1 & -1 & -1 &  1 & -1 & -1 & -1 &  1\\
++--+ &  1 & -1 & -1 & -1 & -1 &  1 &  1 &  1 & -1 & -1 &  1 &  1 & -1 & -1 &  1\\
++-+- & -1 &  1 & -1 & -1 &  1 & -1 &  1 & -1 &  1 & -1 &  1 & -1 &  1 & -1 &  1\\
++-++ &  1 &  1 &  1 & -1 & -1 & -1 &  1 &  1 &  1 & -1 &  1 &  1 &  1 & -1 &  1\\
+++-- & -1 & -1 &  1 &  1 & -1 & -1 &  1 & -1 & -1 &  1 &  1 & -1 & -1 &  1 &  1\\
+++-+ &  1 & -1 & -1 &  1 &  1 & -1 &  1 &  1 & -1 &  1 &  1 &  1 & -1 &  1 &  1\\
++++- & -1 &  1 & -1 &  1 & -1 &  1 &  1 & -1 &  1 &  1 &  1 & -1 &  1 &  1 &  1\\
+++++ &  1 &  1 &  1 &  1 &  1 &  1 &  1 &  1 &  1 &  1 &  1 &  1 &  1 &  1 &  1
\end{array}
$
\end{sideways}
\normalsize
\end{table}
In this case the run of \texttt{hilbert} took 1008.97 seconds, and the number
of elements of the Hilbert basis is 26142. Some of the arrays are
quite unusual. Table \ref{tab:4}
\begin{table}
  \centering
  \caption{Support (supp), total, maximum replication (maxrep), for the Hilbert basis of the OA$(5,2)$}
  \label{tab:4}
  \begin{tabular}{l|rrrrrrrr|rrrr|}
& \multicolumn{8}{c}{total} & \multicolumn{4}{c}{maxrep} \\
& 8 & 12 & 16 & 20 & 24 & 28 & 32 & 36 &   1&   2&   3&   4\\
\hline
supp  \\
8 & 60 & 0 & 0 & 0 & 0 & 0 & 0 & 0 & 60 & 0 & 0 & 0\\
11 & 0 & 32 & 0 & 0 & 0 & 0 & 0 & 0 & 0 & 32 & 0 & 0\\
12 & 0 & 192 & 0 & 0 & 0 & 0 & 0 & 0 & 192 & 0 & 0 & 0\\
13 & 0 & 0 & 0 & 480 & 0 & 0 & 0 & 0  & 0 & 480 & 0 & 0\\
15 & 0 & 0 & 0 & 0 & 1920 & 0 & 0 & 0 & 0 & 0 & 1920 & 0\\
16 & 0 & 0 & 162 & 0 & 0 & 2624 & 5760 & 2880 & 162 & 0 & 8064 & 3200\\
18 & 0 & 0 & 0 & 0 & 5760 & 5760 & 0 & 0 & 0 & 5760 & 5760 & 0\\
19 & 0 & 0 & 0 & 480 & 0 & 0 & 0 & 0 & 0 & 480 & 0 & 0\\
21 & 0 & 0 & 0 & 0 & 0 & 0 & 0 & 32 & 0 & 32 & 0 & 0\\
\hline
maxrep\\
1 & 60 & 192 & 162 & 0 & 0 & 0 & 0 & 0\\
2 & 0 & 32 & 0 & 960 & 5760 & 0 & 0 & 32\\
3 & 0 & 0 & 0 & 0 & 1920 & 8064 & 5760 & 0\\
4 & 0 & 0 & 0 & 0 & 0 & 320 & 0 & 2880\\
  \end{tabular}
\end{table}
shows the distribution of the support, the totals and the maximum number of replication in each element of the Hilbert basis.
Let us consider first the 60+192+162 = 414 elements of the Hilbert basis with no replications. If we compare this table with the results mentioned in Section 1, we see that the cases with 8 or 12 points are the same number, and actually the same fractions. In the case with 16 points, 552 - 162 =  390 are missing. This means that in this case all the extra fractions are actually the union of two fractions with 8 points and disjoint support. To check this, we find that 450 couples of elements of the basis with 8 elements have disjoint support.
\subsection{OA$(6,2)$} This case appears to be not solvable with this technology because of the long computational time. The authors of the software we are using announced the future release of an option to the program \texttt{hilbert} that computes only the 0-1 vectors. In our case, this would compute only the OA's with no replication. The examples above show why we expect this number to be much smaller than the number of elements in the full Hilbert basis and, as a consequence, to be computable in practical times. A different option consist in the search of a basis of the OA$(6,2)$ such that prescribed entries are zero. The assignment to zero of elements could be done by random sampling, or in a systematic way. As an example, if we force ten zeros to the first entries with the ordering of treatments as in the previous examples, then the Hilbert basis has 6 elements.
\section{Discussion}
The case OA$(6,2)$ could not be computed because of the excessive computational time. As many other general methods, this one is currently restricted to small example, and has to be considered mainly of conceptual interest. The main problem is the big number of fractions with replications bigger than 1. A fraction with no replication can be generated only by elements of the basis of the same type, then the algorithm could possibly be modified to compute only the subset generating fractions without replications bigger than 1.

Constrains of the form $\sum_{a \in D} R(a)= T$ for a given total $T$, or $R(a)=1$ for some $a\in D(n)$, change the character of th problem and require different algorithms and software. The second constrain is of particular interest, because reduces the cardinality due to symmetries. In the polynomial representation, this constrain reduces to $\sum_{\alpha \in L} b_\alpha = 1$. Constrains of the form $\sum_{a \in D} R(a) \propto T$ for a given total $T$ do not reduce the computational complexity.

\bibliographystyle{jstp}
\bibliography{/Volumes/Debian/gianni/Archive/bibs/tutto}
\end{document}